\documentclass[11pt,a4paper]{article}
\usepackage[]{times}
\usepackage{fullpage}
\usepackage{graphicx}
\usepackage{float}
\usepackage{subfigure}
\usepackage{amsmath}
\usepackage{fancyhdr}
\usepackage[colorlinks]{hyperref}
\usepackage{xcolor}
\date{}%leave empty

\newcommand{\natu}{{\sf I \! N}}
\newtheorem{e-proposition}[theorem]{Proposition}

\newtheorem{e-definition}[theorem]{Definition\rm}

\title{A comparative study of a higher order finite element method to    
\\     handle natural boundary conditions for smooth spacial domains}
\author{
    Vitoriano Ruas$^{1}$\thanks{Sorbonne Universit\'e, Campus Pierre et Marie Curie, 4 place jussieu, Couloir 55-65, 4\`eme \'etage, 75005 Paris, France.}
		\\[1mm]
  {\small $^{1}$ Institut Jean Le Rond d'Alembert, CNRS UMR 7190, 
  Sorbonne Universit\'e, Paris, France.}\\[1mm]
  {\small e-mail: {\it vitoriano.ruas@upmc.fr}}}

\begin{document}

\maketitle

\begin{abstract}
A numerical study of tetrahedral Lagrange finite element methods upon polynomials of degree greater than one is presented in the solution of second order boundary value problems posed in a curved spatial domain, with Neumann or Robin conditions. In this case the isoparametric technique combined with a standard Galerkin formulation is generally considered to produce a satisfactory accuracy improvement, as compared to the so-called do-nothing strategy based on straight-edged simplicial fitted meshes. In a recent paper the author proposed an alternative technique to the latter method in the two-dimensional case, in which corrective terms are added to the underlying variational form, in such a way that the prescribed natural conditions are enforced on the true curvilinear boundary instead of the one of the approximating polygon. This shift - also known as an off-site enforcement -, gives rise to a Petrov-Galerkin formulation, whose reliability was formally established in that article and also thoroughly certified by comparative numerical experiments. The purpose of this work is two-fold: First of all a rather straightforward extension of the new technique to the three-dimensional case is described, in connection with fitted straight-edged tetrahedral meshes. Then most of the paper is dedicated to numerical comparisons of this approach with both the classical isoparametric technique and the do-nothing strategy, in the case of quadratic elements. 
  

\end{abstract}

\section{Introduction}

Effective ways to approximate the solution of boundary value problems posed in curved domains via discretization methods such as the finite difference or the finite element method was a topic extensively addressed in the literature as early as in the middle of the last century. In the framework of the latter the isoparametric technique has been widely used by practitioners up to nowadays, as far as second order equations are concerned (cf. \cite{Zienkiewicz}, \cite{Ciarlet}, \cite{CiarletRaviart}). However, its extension to solution methods for other types of problems is not straightforward, if not known yet. Maybe for this reason the subject is still deserving the attention of specialists (see e.g. \cite{FBMS}, \cite{FBRT}, \cite{Arxiv2017}, \cite{ZAMM}, \cite{IMAJNA}, \cite{Scovazzi}, \cite{M2AN}, \cite{DGS}) up to now. In recent work, most of the effort was made to handle essential boundary conditions, i.e., the case where degrees of freedom are prescribed on a curvilinear boundary. The application of isoparametric elements to the case of natural boundary conditions such as Neumann's and Robin's in turn was studied by several authors in the seventies and the eighties (see e.g. \cite{BarrettElliott} and references therein). However, to the best of the author's knowledge, since then no other approach had been considered to treat this case until he proposed himself an alternative. This is a Petrov-Galerkin formulation based both on fitted straight-edged simplicial meshes and polynomial algebra, in which the natural conditions are enforced off-site, that is, on the boundary of the true domain. A rather complete theoretical and numerical study of this new technique in the two-dimensional case can be found in \cite{Arxiv2026}. Here, after showing that its extension to the case of straight-edged tetrahedral meshes is straightforward, the resulting formulation was implemented and assessed.\\
Actually, in this article we compare the aforementioned new formulation with both the isoparametric technique and the do-nothing strategy. By the latter denomination we mean the solution method in which the curved domain is approximated by the union of the straight-edged tetrahedra of a fitted mesh and the same Galerkin formulation as for the curved domain is employed in the resulting approximated domain. We carry out our numerical study for tetrahedral Lagrange finite element methods in the solution of second order partial differential equations posed in a smooth three-dimensional domain with both Neumann and Robin homogeneous or inhomogeneous boundary conditions. More specifically, the three approaches to handle the mismatch of the boundary of the domain formed by tetrahedral finite-element meshes and the curved domain are compared. %Two of them are classical methods known since long, while the third one is introduced and implemented for the first time in this work. 
As a matter of fact, while on the one hand such a study does not make much sense in the case of affine finite elements, it becomes relevant as long as a higher order method is employed. Thus, as an initial step to confront the three methods with each other, we observe their performances here in terms of both accuracy and processing time among other criteria in connection with quadratic finite elements only. Moreover, although all the three techniques under consideration apply to different kinds of equations, the study is carried out for a model linear reaction-diffusion problem. \\

An outline of the article is as follows. In Section 2, after specifying the model problem to be solved in our study, we give some preparatory notations and definitions related to the straight-edged finite element meshes to be used in the sequel. Our alternative finite element method to handle Neumann or Robin boundary conditions for this model equation is described in Section 3. In Section 4 we report a selection of results from those obtained by means of thorough numerical experimentation carried out on a comparative basis with the two aforementioned solution techniques. We conclude in Section 5 with a few comments.  

\section{Preliminaries}

\hspace{4mm} In this section we supply relevant preparatory material.

\subsection{The model problem}

\hspace{4mm} The model problem that we consider in this article is the same reaction-diffusion equation with inhomogeneous Robin boundary conditions as in the work by Barrett and Elliott \cite{BarrettElliott}. This choice enables the assessment of our methodology as compared to the  isoparametric method studied by these authors in this paper of theirs.\\ 
Let us assume that $\Omega$ is a smooth spatial domain with boundary $\Gamma$ of the $C^{k-1}$-class for an integer $k >1$  and that  
$f$ and $g$ are given functions respectively in $H^{k-1}(\Omega)$ and $H^{k-1/2}(\Gamma)$. We also assume that $d$ is a function defined in $\bar{\Omega}$ % $\sigma \in W^{k,\infty}(\Omega)$ with $\sigma_M \geq \sigma \geq \sigma_m >0$,   
fulfilling $d \geq d_0 > 0$, $d_0 \in \Re$, and also that $r$ and $\gamma$ are functions defined in $\bar{\Omega}$ and $\Gamma$ such that $r \geq r_0 \geq 0$ and $\gamma \geq \gamma_0 \geq 0$, where $r_0$ and $\gamma_0$ are real numbers satisfying $r_0 + \gamma_0 >0$. Moreover, since our method is basically designed to deal with quite smooth solutions, we shall assume that $\Omega$, $d$, $r$ and $\gamma$ enjoy the same regularity properties as in \cite{BarrettElliott}. \\
In doing so, the problem to solve is   
\begin{equation}
\label{eq}
\left\{
\begin{array}{l}
\mbox{Find } u \in H^{k+1}(\Omega) \mbox{ such that} \\
- \nabla \cdot ( d \nabla u) + r u = f \; \mbox{ in } \Omega \\
\mbox{and } \displaystyle d \frac{\partial u}{\partial n} + \gamma u = g \mbox{ on } \Gamma.
\end{array}
\right. 
\end{equation}
In case $\Omega$ is not convex, at least in theory, it will be necessary to extend the data $f$, $d$ and $r$ to $\tilde{\Omega} \setminus \Omega$, where $\tilde{\Omega}$ is a domain fully containing $\bar{\Omega}$, though not much larger than $\Omega$. The resulting extensions denoted by $\tilde{f}$, $\tilde{d}$ and $\tilde{r}$ are assumed to satisfy $\tilde{f} \in H^{k-1}(\tilde{\Omega})$, $\tilde{d} \in C^k(\bar{\tilde{\Omega}})$ and $\tilde{r} \in C^{k-1}(\bar{\tilde{\Omega}})$. These extensions are aimed at enabling the proper definition of a finite element analog of \eqref{eq} in an approximate domain not fully contained in $\Omega$, as seen hereafter. \\
We assume that there exists a strictly positive constant $\tilde{d}_0$ and a non negative constant $\tilde{r}_0$ such that $\tilde{d} \geq \tilde{d}_0$ and $\tilde{r} \geq \tilde{r}_0$ all over $\tilde{\Omega}$, with $\tilde{r}_0 + \gamma >0$. For convenience we also consider that $u$ is extended to the whole $\tilde{\Omega}$ by a function $\tilde{u} \in H^{k+1}(\tilde{\Omega})$, for example, in the way advocated in \cite{Stein}. Notice, however, that $\tilde{u}$ is by no means expected to solve of a reaction-diffusion equation like \eqref{eq} in $\tilde{\Omega}$ with $d$, $r$ and $f$ replaced with $\tilde{d}$, $\tilde{r}$ and $\tilde{f}$. However, as seen hereafter, the knowledge of $\tilde{f}$ outside $\Omega$ is not indispensable. Thus we may assume that $\tilde{f}:= -\nabla \cdot (\tilde{d}\nabla \tilde{u}) + \tilde{r}\tilde{u}$ in $\tilde{\Omega}\setminus \Omega$, even if $\tilde{u}$ is unknown in $\tilde{\Omega} \setminus \Omega$.     

\subsection{Meshes, related sets and notations}

\hspace{4mm} Next we give some notations related to fitted finite-element meshings of $\Omega$. \\ 
Let ${\mathcal P} := \{{\mathcal T}_h\}_h$ be a family of meshes of $\Omega$ consisting of straight-edged tetrahedra satisfying the usual compatibility conditions for the finite element method. Defining $h_T$ to be the maximum edge length of a given tetrahedron $T \in {\mathcal T}_h$, $h$ represents the maximum of $h_T$ over ${\mathcal T}_h$. The mesh indexing parameter $h$ is defined by 
$$h:= \displaystyle \max_{T \in {\mathcal T}_h} h_T.$$
Every tetrahedron of ${\mathcal T}_h$ is considered to be a closed set and ${\mathcal T}_h$ is assumed to fit $\Omega$ in such a way that all the vertexes of the polyhedron $\Omega_h$ lie on $\Gamma$, where $\Omega_h$ is the interior of $\cup_{T \in {\mathcal T}_h} T$. %We further define $\tilde{\Omega}_h := \Omega \cup \Omega_h$ and $\Omega^{'}_h := \Omega \cap \Omega_h$. 
The boundary of $\Omega_h$ %and $\tilde{\Omega}_h$ are respectively 
is a piecewise triangular surface denoted by $\Gamma_h$, whereas $\tilde{\Omega}_h$ is the set $\Omega \cup \Omega_h$.\\
We must stress the fact that here the family ${\mathcal P}$ of triangulations of $\Omega$ is not assumed to be quasi-uniform. Indeed, akin to \cite{Arxiv2026}, in this article we only require ${\mathcal P}$ to be a regular family of meshes in the sense of \cite{Ciarlet}. Denoting the shortest edge of $T$ by $l_T$, this means that there exists a strictly positive constant $\mu$ such that $l_T \geq \mu h_T$, $\forall T \in {\mathcal T}_h$ and $\forall {\mathcal T}_h \in {\mathcal P}$. \\

In this work the notation $|\cdot|$ is employed to represent either the volume of a subset of $\Re^3$ with a non empty interior or the area of a plane surface or yet the length of a straight segment contained in $\tilde{\Omega}_h$.\\     
Assuming that every element in ${\mathcal T}_h$ has at most one face contained in $\Gamma_h$, let ${\mathcal S}_{h}$ be the subset of ${\mathcal T}_h$ consisting of tetrahedra $T$ having one face on $\Gamma_{h}$, say $F_T$. Referring to Figure 1, for every $T \in {\mathcal S}_{h}$ we denote by $O_T$ the vertex of $T$ not belonging to $\Gamma$ and by ${\bf n}_T$ the unit outer normal vector to $\Gamma_{h}$.\\ 
Now, for every $M \in F_T$, let $P$ be the intersection with $\Gamma$ in the immediate neighborhood of $T$ of the line perpendicular to $F_T$ passing through $M$. We assume that for all meshes in ${\mathcal P}$, $h$ is sufficiently small for $P$ such that $|\overline{MP}| \leq h_T$ to be unique $\forall M \in F_T$ and $\forall T \in {\mathcal S}_h$. We further denote by ${\bf n}_P$ the unit outer normal vector to $\Gamma$ at $P$, as illustrated in Figure 1. Under the above assumption on $h$, let $\Gamma_T$ be the set of all such points $P$ for a given $T \in {\mathcal S}_h$. 
\begin{figure}[h!]
\label{fig1}
\centerline{\includegraphics[width=3.8in]{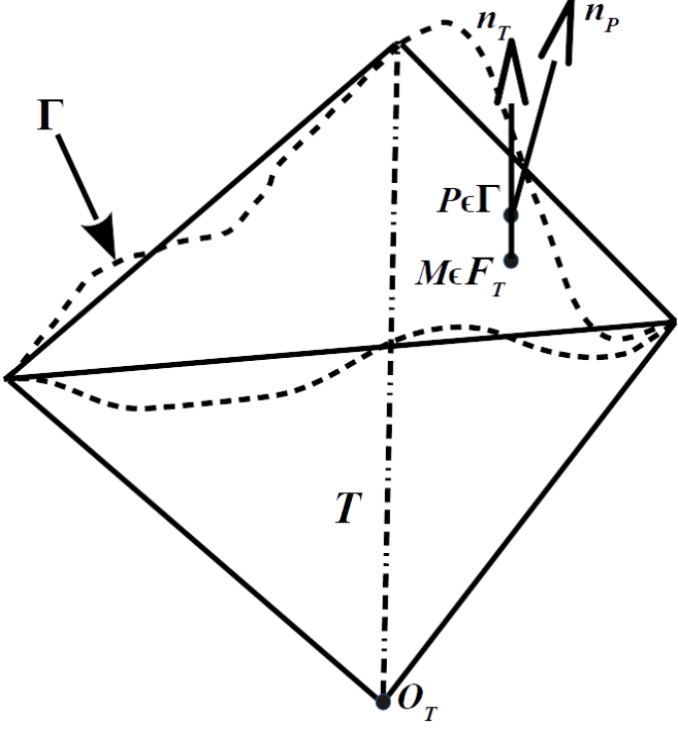}}
\vspace*{8pt}
\caption{A point $M \in F_T$ and normal vectors ${\bf n}_T$ at $M$ and ${\bf n}_P$ at $P \in \Gamma$ with $\overline{MP} \perp F_T$ for $T \in {\mathcal S}_{h}$}
\end{figure} 
and further define $\Delta_T$ to be the subset of $\tilde{\Omega}_h$ delimited by $F_T$, $\Gamma_T$ and three subsets of the planes orthogonal to $F_T$ containing the edges of this face delimited by $\Gamma_T$ and the respective edge. We also set $\tilde{T} := T \cup \Delta_T$ $\forall T \in {\mathcal S}_h$ and $\tilde{T} = T$ $\forall T \in {\mathcal T}_h \setminus {\mathcal S}_h$.\\
%As for additional notations related to the mesh ${\mathcal T}_h$, the polygonal domain $\Omega_h$ 
%we extend $\| \cdot \|_{1,h}$ to the standard norm of either $H^j(\Omega_h)$ for $j > 1$ or $L^2(\Omega_h)$ for $j=0$. Moreover, in order to accommodate non convex domains, we introduce the set $\tilde{\Omega}_h := \Omega \cup \Omega_h$. We also represent by $\tilde{T}$ the set $T \cup \Delta_T$ $\forall T \in {\mathcal S}_h$ and define $\tilde{T} := T \; \forall T \in {\mathcal T}_h \setminus {\mathcal S}_h$.
 %Further $D^j{\mathcal A}$ represents the $j$-th order tensor, whose components are the $j$-th order partial derivatives with respect to the space variables of a function ${\mathcal A}$ in the strong or in the weak sense.
Finally, $\theta_T(M)$ is defined to be the angle between ${\bf n}_T$ and ${\bf n}_P$ $\forall M \in F_T$ and $\forall T \in {\mathcal S}_h$. Since $h$ is supposedly small, we may also assume that ${\mathcal P}$ is such that there exist two mesh-independent constants $s_m$ and $c_m$ fulfilling 
\begin{equation}
\label{Thetam}
\left\{
\begin{array}{l}
\max_{T \in {\mathcal S}_h} \max_{M \in F_T} | \theta_T(M)| \;  \leq \; s_m h \implies 
\max_{T \in {\mathcal S}_h} \max_{M \in F_T} | sin \theta_T(M)| \; \leq \; s_m h \\  
\mbox{and } \min_{T \in {\mathcal S}_h} \min_{M \in F_T} | cos \theta_T(M)| \; \geq \; c_m > 0 \; \forall {\mathcal T}_h \in {\mathcal P}.
\end{array}
\right.
\end{equation} 
In complement to the above definitions, assumptions and notations we shall use the following ones. \\
$P_k(D)$ is defined to be the space of polynomials of degree less than or equal to $k$ restricted to a subset $D \in \Re^n$, for $n \in \natu$.\\
Hereafter we denote by $\| \cdot \|_{0,h}$ the standard norm of $L^2(\Omega_h)$ and by $\| \cdot \|_{1,h}$ and ]$| \cdot |_{1,h}$ the standard norm and semi-norm of $H^1(\Omega_h)$. Analogously, we employ the notations $\| \cdot \|_{0,\tilde{h}}$ and $| \cdot |_{1,\tilde{h}}$ for the respective norms  calculated by replacing $\Omega_h$ with the curved domain approximating $\Omega$ associated with isoparametric finite elements.

\section{A new approach to handle natural boundary conditions}

In this section we go into the description of the Petrov-Galerkin formulation advocated in this work to solve \eqref{eq}, which yields a finite-element analog of this equation.\\
\indent First of all, for an integer $k >1$, we denote by $V_h^k$ the space of standard Lagrange finite elements associated with ${\mathcal T}_h$, i.e., the space of continuous functions, whose restriction to every $T \in {\mathcal T}_h$ is a polynomial of degree less than or equal to $k$. \\
Let $\Delta_{\Omega}$ be the subset of $\Omega$ defined by 
$\Delta_{\Omega}: = \displaystyle \cup_{T \in {\mathcal S}_h} \Delta_T \setminus \Omega_h.$  
We extend every $v \in V^k_h$ to $\Delta_{\Omega}$ as follows: $\forall T \in {\mathcal S}_h$ the expression of $v$ in $T$ in terms of the space variables is applied as such to all points in $\Delta_T$ not belonging to $T$. \\
We represent by $\tilde{V}^k_h$ the space of functions $\tilde{v}$ defined in $\tilde{\Omega}_h$ which are continuous in $\Omega_h$ such that $\forall T \in {\mathcal T}_h$ the restriction of $\tilde{v}$ to $\tilde{T}$ belongs to $P_k(\tilde{T})$.\\
Now with every function $\varphi$ defined on $\Gamma$ we associate a function $\bar{\varphi}$ uniquely defined on $\Gamma_h$ by $\bar{\varphi}(M) = \varphi(P)$ $\forall M \in F_T$ and $\forall T \in {\mathcal S}_h$, where $P$ is the point of $\Gamma$ associated with $M \in F_T$ in the way specified in Subsection 2.2.      
Recalling the unit vectors ${\bf n}_T$ and ${\bf n}_P$ for $M \in F_T$ and $T \in {\mathcal S}_h$, we consider an arbitrarily chosen unit vector ${\bf s}_P$ of the tangent plane to $\Gamma_T$ at $P$ and define a second unit vector of this plane by ${\bf t}_P := {\bf n}_P \times {\bf s}_P$. Let us denote respectively by $\partial \cdot /\partial n_P$, $\partial \cdot /\partial s_P$ and $\partial \cdot /\partial t_P$  the first order partial derivatives in the direction of ${\bf n}_P$, ${\bf s}_P$ and ${\bf t}_P$ and by $\theta_T$, $\sigma_T$ and $\tau_T$ the angles between ${\bf n}_P$, ${\bf s}_P$, ${\bf t}_P$ and ${\bf n}_T$. In connection with the these three angles we define the  
functions $\nu_T$, $\lambda_T$ and $\mu_T$ on $F_T$, given by $\nu_T =$ cos $\theta_T$, $\lambda_T =$ cos $\sigma_T$ and $\mu_T =$ cos $\tau_T$, together with corresponding three functions $\nu$, $\lambda$ and $\mu$ defined on $\Gamma_h$ such that $\nu(M)=\nu_T(M)$, $\lambda(M)=\lambda_T(M)$ and $\mu(M)=\mu_T(M)$ $\forall M \in F_T$ and $\forall T \in {\mathcal S}_h$.\\
Now for all functions $v \in V_h^k$ and $w \in H^3(\tilde{\Omega}_h) + \tilde{V}_h^k$ where, for a natural number $m$ and a bounded open set $D$ of $\Re^3$, $H^m(D)$ is the classical Hilbert-Sobolev space \cite{Adams}, we define a bilinear form $\hat{a}_{h}$ and 
a linear form $\hat{\mathcal F}_{h}$ as follows.
\begin{equation}
\label{hatah}
\left\{
\begin{array}{l}
\hat{a}_h(w,v) := \displaystyle \int_{\Omega_h} (\tilde{d} \nabla w \cdot \nabla v \;+ \; \tilde{r} wv) \\   
+ \displaystyle \sum_{T \in {\mathcal S}_h} \int_{F_T} \left\{\nu_T \left[ \bar{\gamma} \bar{w} \; + \; \displaystyle \int_M^P \frac{\partial}{\partial n_T}\left(\tilde{d} \frac{\partial w}{\partial n_P}\right)\right] \; + \; \tilde{d} \left( \lambda_T \frac{\partial w}{\partial s_P} + \mu_T \frac{\partial w}{\partial t_P} \right) \right\}v
\end{array}
\right.
\end{equation}
\begin{equation}
\label{hatFh}
\hat{\mathcal F}_h(v):=\displaystyle \int_{\Omega_h} \tilde{f} v \; + \; \displaystyle \sum_{T \in {\mathcal S}_h} \int_{F_T} \nu_T \bar{g} v.
\end{equation}
A natural way to define our finite element approximation of \eqref{eq} is to 
solve the following problem
\begin{equation}
\label{hatvph}
\mbox{Find } \hat{u}^k_h \in \tilde{V}^k_h \mbox{ such that } \hat{a}_{h}(\hat{u}^k_h,v) = \hat{\mathcal F}_h(v) \; \forall v \in V^k_h. 
\end{equation}
The main asset of the discrete variational problem \eqref{hatvph} is the fact that, at least in the case where $\Omega$ is convex, the underlying residual $\hat{a}_h(u,v)-\hat{\mathcal F}_h(v)$ vanishes $\forall v \in V_h^k$. This property can be established in the same manner as in \cite{Arxiv2026} for the two-dimensional case. Actually it also holds in the case of non convex domains, provided $u$ is replaced by its extension $\tilde{u}$ to the slightly larger domain $\tilde{\Omega}$ considered in Subsection 2.1.\\  

As we should note, however, evaluating both $\hat{a}_h$ and $\hat{\mathcal F}_h$ exactly may be unpractical in some cases if not unfeasible. Therefore we will pose instead a similar variational formulation in which the data $\tilde{d}$, $\tilde{r}$, $\tilde{f}$ restricted to every $T$ in ${\mathcal T}_h$ and the data $\nu$, $\lambda$, $\mu$, $\bar{\gamma}$ and $\bar{g}$, together with $\overline{\hat{u}}_h^k$, all  restricted to $F_T$ for every $T$ in ${\mathcal S}_h$ are replaced by suitable Lagrange interpolation polynomials $q_{\tilde{d}}$, $q_{\tilde{r}}$, $q_{\tilde{f}}$ and $s_{\nu}$, $s_{\lambda}$, $s_{\mu}$, $s_{\bar{g}}$, $s_{\bar{\gamma}}$ and $s_{\bar{u}}$. The guiding principle for the choice of such interpolating polynomials is two-fold. On the one hand the deviation from the exact solution brought about by this additional error source must be compatible with the expected order of convergence of the formulation \eqref{hatvph}. On the other hand, the interpolation points in the tetrahedra must lie in their interior sufficiently away from $\Gamma_h$, in order to make sure that the functions to interpolate are well defined at all of them. Rigorously, the latter restriction applies only to $f$ and to tetrahedra $T \in {\mathcal S}_h$ such that $T \setminus \Omega \neq \emptyset$. However, for coherence and also for the sake of simplicity, we choose the same interpolation procedure for all tetrahedra in the mesh. \\
In doing so, we only have to integrate polynomials in $T$, on $F_T$ and along the chords $\overline{MP}$. Therefore, strictly speaking, we can do without numerical quadrature, in contrast to the case of isoparametric elements.\\
Nevertheless, akin to \cite{Arxiv2026} we will not use polynomial interpolation of the data in this work, but instead numerical integration  in the tetrahedra, the boundary faces and the chords. This is because our method's implementation becomes fairly simpler in this case, even if we admit that, rigorously, the regularity requirements for the data in order to attain optimal orders of convergence become much more severe for $k >3$ (cf. \cite{Arxiv2026}).\\    
More concretely, for $\bar{a}_h$, in every $T \in {\mathcal T}_h$ we use a $J_k$-point Gauss quadrature formula with points $Q_j \in T$ and strictly positive weights $\chi_j$ for $j=1,\ldots J_k$ for non polynomial integrands in $T$, an $I_k$-point Gaussian quadrature formula for integrals on $F_T$ with points $M_i$ and strictly positive weights $\omega_i$ for $i=1,\ldots I_k$ and an $L_k$-point Gaussian quadrature formula for integrals along the chords $\overline{M_iP}$ with points $N_l$ and strictly positive weights $\xi_l$ for $l=1,\ldots L_k$. As usual, we assume that the sums of the weights $\chi_j$, $\omega_i$ and $\xi_l$ are all equal to one. Besides the latter property we have to be more specific about these quadrature formulae. For the moment we just state the following:\\
The integrals of a given continuous function $\phi$ in $T \in {\mathcal T}_h$ is approximated by ${\mathcal J}_T^k(\phi)$, where 
\begin{equation}
\label{JTk} 
{\mathcal J}_T^k(\phi) := |T| \displaystyle \sum_{j=1}^{J_k} \chi_j \phi(Q_j) 
\end{equation}
In order to do without the explicit knowledge outside of $\Omega$ of the extensions $\tilde{f}$ and $\tilde{r}$ of the data $f$ and $r$, we assume that for all $T \in {\mathcal S}_h$ all the quadrature points $Q_j$ lie inside $\Omega$, in case this domain is not convex.\\
The integral of $\phi$ on $F_T$ for $T \in {\mathcal S}_h$ in turn are approximated by ${\mathcal I}_T^k(\phi)$, where           
\begin{equation}
\label{ITk} 
{\mathcal I}_T^k(\phi) := |F_T| \displaystyle \sum_{i=1}^{I_k} \omega_i \phi(M_i). 
\end{equation} 
Finally, the integral of $\phi$ along the chords $\overline{M_iP}$ where $M_i$ is a quadrature point on $F_T$ in the formula \eqref{ITk} for $T \in {\mathcal S}_h$ are approximated by ${\mathcal L}_{T,i}^k(\phi)$, where           
\begin{equation}
\label{LTk} 
{\mathcal L}_{T,i}^k(\phi) := |\overline{M_iP}| \displaystyle \sum_{i=1}^{L_k} \xi_l \phi(N_l). 
\end{equation}  
Notice that, in view of the definition \eqref{hatah}, if $\Omega$ is not convex we need the explicit knowledge of an extension $\tilde{d}$ of $d$ outside $\Omega$ in order to apply the formula \eqref{ITk} for some $T \in {\mathcal S}_h$ and formula \eqref{LTk} for those chords not completely immersed in $\Omega$. \\

Now for every $w \in H^3(\tilde{\Omega}_h) + \tilde{V}_h^k$ and $v \in V_h^k$ we set
\begin{equation}
\label{ah}
%\left\{
%\begin{array}{l}
a_h(w,v) := b_h(w,v) + c_h(w,v), 
%\end{array}
%\right. 
\end{equation}
where $\forall (w;v) \in [H^3(\tilde{\Omega}_h) + \tilde{V}^k_h] \times V^k_h$
\begin{equation}
\label{bh}
b_h(w,v) := \; \displaystyle \sum_{T \in {\mathcal T}_h} {\mathcal J}_T^k \left( \tilde{d} \nabla w \cdot \nabla v \; + \; \tilde{r} w v \right) 
%\; + \; \displaystyle \sum_{T \in {\mathcal S}_h} {\mathcal I}_T^k \left(\nu_T \bar{\gamma} \bar{w} v \right)  
\end{equation}
and
\begin{equation}
\label{ch}
%\left\{
%\begin{array}{l}
c_h(w,v) := \; \displaystyle  \sum_{T \in {\mathcal S}_h} {\mathcal I}_T^k  \left( \left[ \nu_T \left\{ \bar{\gamma} \bar{w} \; + \; \displaystyle {\mathcal L}_{T,i}^k \left(\frac{\partial}{\partial n_T} \left[\tilde{d} \frac{\partial w}{\partial n_P}\right] \right) \right\} \; + \; \tilde{d} \displaystyle \left\{ \lambda_T \frac{\partial w}{\partial \sigma_P} + \mu_T \frac{\partial w}{\partial \tau_P} \right\} \right] v \right). 
%\right. 
\end{equation}
As far as $\hat{\mathcal F}_h$ is concerned, we define its approximation ${\mathcal F}_h$ by
\begin{equation}
\label{Fh}
{\mathcal F}_{h}(v) := \; \displaystyle \sum_{T \in {\mathcal T}_h} {\mathcal J}_T^k(f v) + \; \displaystyle \sum_{T \in {\mathcal S}_h}
{\mathcal I}_T^k \left(\nu_T [\bar{g}]_{|F_T} v \right) \; \forall v \in V^k_h.
\end{equation}

With the above definitions, we set the approximate problem to solve as
\begin{equation}
\label{vph}
\mbox{Find } u^k_{h} \in \tilde{V}^k_h \mbox{ such that } a_{h}(u^k_{h},v) = {\mathcal F}_{h}(v) \; \forall v \in V_h^k,
\end{equation}
Notice that both $V^k_h$ and $\tilde{V}^k_h$ are finite-dimensional spaces with the same dimension, Hence both $a_h$ and ${\mathcal F}_h$ are continuous forms over $\tilde{V}^k_h \times V^k_h$ and $V^k_h$ equipped with any norm. For example, we have
\begin{equation}
\label{continuityah}
a_h(w,v) \leq A_h \| w \|_{1,h} \| v \|_{1,h} \; \forall (w;v) \in \tilde{V}^k_h \times V^k_h
\end{equation}
and
\begin{equation}
\label{continuityFh}
{\mathcal F}_h(v) \leq C_h \| v \|_{1,h} \; \forall v \in V^k_h,
\end{equation}
for suitable constants $A_h$ and  $C_h$.\\
According to \cite{Babuska1}, it follows that problem \eqref{vph} is well posed if the bilinear form $a_h$ satisfies an $inf-sup$ condition with a constant $\alpha > 0$ over $\tilde{V}^k_h$ and $V^k_h$ equipped with any norm. In the case under study the natural norm is $\| \cdot \|_{1,h}$, for which this condition writes  
\begin{equation}
\label{infsuph}
\forall w \in \tilde{V}^k_h \displaystyle \sup_{v \in V^k_h \setminus \{0\}} \frac{a_h(w,v)}{\| v \|_{1,h}} \geq \alpha \| w \|_{1,h}.
\end{equation}      
Notice that conditions similar to \eqref{continuityah} and \eqref{continuityFh} hold for $\hat{a}_h$ and $\hat{\mathcal F}_h$. \\
The proof of condition \eqref{infsuph} in the analogous two-dimensional case given in \cite{Arxiv2026}, with $\alpha$ independent of $h$ as required, is very laborious. In this work we adopt a different approach by assuming that problem \eqref{vph} is well-posed and focus on checking out in several test-cases that such an assumption holds true. Moreover we also observe that optimal orders of convergence are attained at least for $k=2$, by means of a judicious choice of the numerical quadrature formulae.   \\

To close this section we recall the problem to solve for the so-called do-nothing strategy. 
Let $\bar{a}_{h}$ be the bilinear form $b_h$ defined by \eqref{bh} on $V^k_h \times V^k_h$. We further define 
\begin{equation}
\label{barFh}
\bar{\mathcal F}_{h}(v) := \; \displaystyle \sum_{T \in {\mathcal T}_h} {\mathcal J}_T^k(f v) + \; \displaystyle \sum_{T \in {\mathcal S}_h}
{\mathcal I}_T^k\left([\bar{g}]_{|F_T} v \right) \; \forall v \in V^k_h.
\end{equation} 
Then we solve the problem 
\begin{equation}
\label{barvph}
\left\{
\begin{array}{l}
\mbox{Find } \bar{u}_{h}^k \in V^k_h \mbox{ such that}\\
\bar{a}_h(\bar{u}^k_{h},v) = \bar{\mathcal F}_h(v) \; \forall v \in V^k_h.
\end{array}
\right.
\end{equation}
   
\section{Numerical experiments}

In this section we go into numerical experimentation, which is the main purpose of this paper. All the test-problems considered below were  solved by three methods associated with quadratic elements ($k=2$), namely, the Petrov-Galerkin formulation \eqref{vph}, the standard isoparametric technique and the do-nothing strategy \eqref{barvph}. \\

\subsection{Preliminaries} 
We checked out the performances of the three methods by solving test-problems with constant data $d$, $r$ and $\gamma$ and a known exact solution $u$ satisfying three different types of boundary conditions, namely, homogeneous Neumann and inhomogeneous Neumann and Robin boundary conditions. The data $f$ and $g$ are determined from \eqref{eq} accordingly.\\
In all cases $\Omega$ is a domain equal to or strictly contained in the ellipsoid $\Omega_{ab}$ with semi-axes equal to $a$, $b$ and $1$ with $a \leq 1$ and $b \leq 1$. More concretely, we consider both the case where $\Omega$ is $\Omega_{ab}$ itself or the case where $\Omega$ is the non convex domain $\Omega_{\beta}$ defined by a real parameter $\beta \in (0,1/2]$ to be a perturbation of the unit ball $\Omega_{11}$ centered at the origin of the coordinate axes. More precisely, the boundary of $\Omega_{\beta}$ is expressed in spherical coordinates $(\rho,\varphi,\theta)$ by $\rho  = 1 - \beta \mbox{ cos }4\varphi$.\\
Only test-problems with three planes of symmetry are considered here. For this reason we confined ourselves to the generation of meshes for the sole octant of $\Omega$ given by $x>0$, $y>0$ and $z>0$. For this purpose we used the procedure described in \cite{LBRSR} defined by an even integer $L$ consisting of $6L^3$ tetrahedra, thereby yielding a quasi-uniform family of meshes parametrized by $L$. Here we took $L = 2^m$ for integers $m$ ranging between $1$ and $4$ in all cases. For better readability we redefine the mesh parameter $h$ as the quantity $h:=1/L$.\\

Hereafter the notation $| v |_{\infty,h}$ stands for the maximum of the values of a function $v$ in $V_h^2 + H^2(\Omega_h)$ at the nodes  of the mesh corresponding to the standard degrees of freedom of $V_h^2$ (cf. \cite{Ciarlet}). Analogously, we employ the   
notation $| v |_{\infty,\tilde{h}}$ for the semi-norm defined upon the nodes associated with the isoparametric technique.\\

Fifth order numerical quadrature formulae ${\mathcal J}_T^2$ and ${\mathcal I}_T^2$ were used to integrate non polynomial functions, namely, a fifteen-point formula for the tetrahedra $T$ (cf. \cite{Stroud, Zienkiewicz}) and a seven-point formula for the triangular faces $F_T$ (see e.g. \cite{Dunavant}). Notice that no quadrature formula ${\mathcal L}_{T,i}^2$ is necessary here, since the integrands along the chords $\overline{M_iP}$ are constant in all cases. \\
  
As for the practical solution of problem \eqref{vph}, we use an auxiliary symmetric coercive Galerkin variant thereof with the same bilinear form $\bar{a}_{h} \equiv b_h$, though with an adjusted right hand side functional at every step of an iterative solution procedure of the same type advocated in \cite{JCAMbis}. The problem to solve at every step leads to a symmetric linear system of algebraic equations with a fixed positive definite matrix. The step-by-step solution procedure functions as follows:\\
Taking the solution $\bar{u}_{h}$ of \eqref{barvph} as the initial guess $u^0_h$ for the solution $u_h$ of \eqref{vph}, with $k=2$ in both cases, for $n=1,2,\ldots$ determine approximations $u_h^n \in V^2_h$ of $u_h$, until a given tolerance $tol$ is satisfied in a sense to be specified, say, for $n=N_{it}$, where $u_h^n$ is the unique solution of the following problem:
\begin{equation}
\label{iterproc}
\bar{a}_{h}(u^n_h,v) = {\mathcal F}_h(v) - c_h(u^{n-1}_h,v) \; \forall v \in V^2_h.
\end{equation}
The above symmetric coercive problem to determine $u^n_h$ for every $n$ can be solved by the classical conjugate gradient method, taking $u^{n-1}_h$ as the initial guess. Since the matrix associated with the problem at hand is very sparse and remains fixed along the iterations, this allows for a reduction of storage requirements to a minimum in the solution of the non symmetric problem \eqref{vph}. Actually, the above solution procedure converges in only a few iterations, thereby bringing about very drastic reductions in terms of running time, as compared to the solution by a direct method, such as Gaussian elimination. \\
For the same reasons as above, the symmetric linear systems with a positive definite matrix to determine $\bar{u}_{h}$ and the isoparametric solution $\tilde{u}_h$ were both solved by the conjugate gradient method starting from a zero initial solution.\\
     
For all the test-problems addressed in this section the numerical results are summarized in tables with the following layout:\\
First we point out that a given table presents results obtained with either the Petrov-Galerkin formulation \eqref{vph} or its Galerkin counterparts \eqref{barvph} with affine mappings and the one with isoparametric (quadratic) mappings. \\
Each table consists of three subsequent groups of three rows. Each one of these groups is associated with the error expressed in one of the three measures $\| \cdot \|_{0,h}$, $| \cdot |_{1,h}$ and $| \cdot |_{\infty,h}$ for \eqref{vph} and \eqref{barvph} and $\| \cdot \|_{0,\tilde{h}}$, $| \cdot |_{1,\tilde{h}}$ and $| \cdot |_{\infty,\tilde{h}}$ for the classical isoparametric solution method.
Each row within a group successively displays the errors for the decreasing values of $h$ obtained by means of the three methods under consideration to approximate the solution $u$ of \eqref{eq} (eventually extended by $\tilde{u}$ outside $\Omega_{\beta}$ in the obvious manner). More specifically, here the resulting approximations are denoted as follows.
\begin{itemize}
\item
$u_h$ standing for the solution of \eqref{vph} with $k=2$;
\item 
$\tilde{u}_h$ for the solution computed with quadratic isoparametric finite elements;
\item
$\bar{u}_{h}$ standing for the solution of \eqref{barvph} with $k=2$.
\end{itemize}
On the right of each error an indicator of the order of convergence ($IOC$) is supplied for the corresponding measure expressed by 
$-\mbox{log}_2 q$ where $q$ is the quotient between the error for the mesh parameter $h$ provided in the corresponding column and the error for the mesh parameter $2h$. Notice that the errors for $2h=1$, i.e., for a mesh of the octant domain with only six tetrahedra, are not supplied. \\

Before going into the numerical results themselves we complete this introduction by giving a few important explanations and definitions, namely,
\begin{itemize}
\item
A stop criterion with a tolerance of $10^{-7}$ for the $| \cdot |_{\infty,h}$ seminorm is applied to both the conjugate gradient algorithm and the iterative procedure to determine $u_h$.
\item
In the tables the $IOC$'s given by $-\mbox{log}_2 q$, where $q$ is the aforementioned quotient between the errors in two successively columns, are denoted respectively by $r_p$, $\tilde{r}_p$ and $\bar{r}_p$ for $u_h$, $\tilde{u}_h$ and $\bar{u}_h$, with $p=0$, $p=1$ and $p=\infty$ in self-explanatory accordance with the corresponding measure. 
\item
In the test-problems whose  exact solution is $u = \rho^2$ with $\rho = \sqrt{x^2+y 2+z^2}$ the $IOC$'s are not supplied for $u_h$ since problem \eqref{vph} yields the exact solution up to the tolerance and round-off errors.
\item
In the sequel we use the following notations for the error functions: $e_h = \tilde{u} - u_h$, $\tilde{e}_h = \tilde{u} - \tilde{u}_h$ and $\bar{e}_h = \tilde{u} - \bar{u}_h$, noting that $\tilde{u} \equiv u$ in the case of $e_h$ and $\bar{e}_h$ for convex domains.
  
\end{itemize}

\subsection{Numerical results}

Next we report a representative selection of the numerical results we have obtained so far in the solution of test-problems of the type \eqref{eq} with $d \equiv 1$, $r \equiv 1$ or $\equiv 0$ and $\gamma \equiv 0$ or $\equiv 1$.\\
The aim of the first three test-problems is to evaluate the methods being compared from the point of view of consistency, by taking a quadratic exact solution.  

\begin{itemize}
\item
\underline{\bf Test-problem 1:} $\Omega=\Omega_{ab}$ with $a=b=1$, $r \equiv 1$, $\gamma \equiv 0$ and $u=\rho^2$ . 

\indent To begin with we check the methods under study in this section for a problem with Neumann boundary conditions. The corresponding results are summarized in Table 1. \\
First of all we observe that the Petrov-Galerkin formulation \eqref{vph} for $k=2$ yields the exact solution as expected, though up to approximation errors of a magnitude close to $10^{-5}$ certainly due to the conjugate effect of the iterative solution procedure and the linear system solver. On the other hand, we recall that the isoparametric method does not enjoy such a property, since it is not based upon piecewise quadratic shape functions. Nevertheless it is very accurate, in particular, as far as the $L^2$-norm is concerned. Moreover 
the $IOC$'s for the three coarsest meshes point to super-convergence with respect to the expected estimates (see e.g. \cite{CiarletRaviart}) of almost half a point in both the norm $\| \cdot \|_{0,h}$ and the semi-norm $| \cdot |_{1,h}$ in this particular case. However, in contrast, this indicator downgrades significantly for the finest mesh. This effect could be due to a tolerance too large for the conjugate gradient solver, as compared to the method's accuracy, although the latter is roughly two orders of magnitude lower in terms of powers of $10$. But, would this be the case, a similar behavior should also be expected of the do-nothing solution $\bar{u}_h$, whose $IOC$'s undoubtedly tend to $2$ for all the three norms under consideration, a value which is the best we can hope for such a simple approach. It should also be noted that, while the accuracy of the latter method is substantially poorer than the one of the other two in terms of the measures $\| \cdot\|_{0,h}$ and $| \cdot |_{\infty,h}$, it is not that far from them as far as the semi-norm $|\cdot|_{1,h}$ is concerned. \\
In short, from the aforementioned observation of the $IOC$ downgrade, we may conjecture that the isoparametric approximation is more quality demanding in terms of the underlying linear system solution, if we wish to exploit in full its theoretical accuracy, especially for finer meshes. We will go back to this discussion within the test-cases that follow.

\begin{table}[h!]
{\small 
\centering
\begin{tabular}{cccccc} &\\ [-.3cm]  
$h$ & $\rightarrow$ & $1/2$ & $1/4$ & $1/8$ & $1/16$  
\tabularnewline &\\ [-.3cm] \hline &\\ [-.3cm]
$\parallel e_h \parallel_{0,h}$ \& $r_0$ & $\rightarrow$ & 0.29807E-6 \& $--$ & 0.49907E-5 \& $--$ & 0.13671E-4 \& $--$ & 0.35368E-5 \& $--$  
\tabularnewline &\\ [-.3cm] \hline &\\ [-.3cm] 
$\parallel \tilde{e}_h \parallel_{0,\tilde{h}}$ \& $\tilde{r}_0$ & $\rightarrow$ & 0.87457E-3 \& $3.61$ & 0.77039E-4 \& $3.50$ & 0.69223E-5 \& $3.48$ & 0.13158E-5 \& $2.40$ 
\tabularnewline &\\ [-.3cm] \hline &\\ [-.3cm] 
$\parallel \bar{e}_h \parallel_{0,h}$ \& $\bar{r}_{0}$ & $\rightarrow$ & 0.12965E+0\& $1.92$ & 0.32818E-1 \& $1.98$ & 0.82296E-2 \& $2.00$ & 0.20590E-2 \& $2.00$ 
\tabularnewline &\\ [-.3cm] \hline &\\ [-.3cm] 
\tabularnewline &\\ [-.3cm] \hline &\\ [-.3cm]
$| e_h |_{1,h}$ \& $r_{1}$& $\rightarrow$ & 0.28505E-6 \& $--$ & 0.18915E-5 \& $--$ & 0.49426E-5 \& $--$ & 0.16681E-4 \& $--$  
\tabularnewline &\\ [-.3cm] \hline &\\ [-.3cm] 
$| \tilde{e}_{h}|_{1,\tilde{h}}$ \& $\tilde{r}_1$ & $\rightarrow$ & 0.12107E-1 \& $1.91$ & 0.24433E-2 \& $2.31$ & 0.46035E-3 \& $2.41$ & 0.85410E-4 \& $2.43$ 
\tabularnewline &\\ [-.3cm] \hline &\\ [-.3cm] 
$| \bar{e}_{h} |_{1,h}$ \& $\bar{r}_1$ & $\rightarrow$ & 0.31900E-1 \& $1.84$ & 0.81970E-2 \& $1.96$ & 0.20633E-2 \& $1.99$ & 0.51695E-3 \& $2.00$ 
\tabularnewline &\\ [-.3cm] \hline &\\ [-.3cm]
\tabularnewline &\\ [-.3cm] \hline &\\ [-.3cm] 
$| e_h |_{\infty,h}$ \& $r_{\infty}$ & $\rightarrow$ & 0.64330E-6 \& $--$  & 0.85352E-5 \& $--$ & 0.21417E-4 \& $--$ & 0.12669E-4 \& $--$ 
\tabularnewline &\\ [-.3cm] \hline &\\ [-.3cm] 
$| \tilde{e}_{h} |_{\infty,\tilde{h}}$ \& $\tilde{r}_{\infty}$ & $\rightarrow$ & 0.59134E-2 \& $3.07$ & 0.69828E-3 \& $3.08$ & 0.79741E-4 \& $3.13$ & 0.14399E-4 \& $2.47$
\tabularnewline &\\ [-.3cm] \hline &\\ [-.3cm] 
$| \bar{e}_{h} |_{\infty,h}$ \& $\bar{r}_{\infty}$ & $\rightarrow$ & 0.19787E+0\& $2.07$ & 0.48785E-1 \& $2.02$ & 0.12154E-1 \& $2.01$ & 0.30377E-2 \& $2.00$ 
\tabularnewline &\\ [-.3cm] \hline &\\ [-.3cm] 
\end{tabular}
%\vskip.2cm 
\caption{Errors \& $IOC$'s for Test-problem 1 solved by \eqref{vph}, isoparametrics and \eqref{barvph} with $k=2$} 
}
\label{table1}
\end{table}
            
\item
\underline{\bf Test-problem 2:} $\Omega=\Omega_{ab}$ with $a=0.6$ and $b=0.8$, $r \equiv 1$, $\gamma \equiv 1$ and $u=\rho^2$ . 

\indent In this test-problem we check the consistency of the three methods in the presence of Robin boundary conditions. \\ 
Essentially the same behavior as for Test-problem 1 can be observed for all of them, as one can infer from Table 2. \\

\begin{table}[h!]
{\small 
\centering
\begin{tabular}{cccccc} &\\ [-.3cm]  
$h$ & $\rightarrow$ & $1/2$ & $1/4$ & $1/8$ & $1/16$  
\tabularnewline &\\ [-.3cm] \hline &\\ [-.3cm]
$\parallel e_h \parallel_{0,h}$ \& $r_0$ & $\rightarrow$ & 0.14182E-6 \& $--$ & 0.12227E-5 \& $--$ & 0.47748E-5 \& $--$ & 0.34884E-5 \& $--$ 
\tabularnewline &\\ [-.3cm] \hline &\\ [-.3cm] 
$\parallel \tilde{e}_h \parallel_{0,\tilde{h}}$ \& $\tilde{r}_0$ & $\rightarrow$ & 0.40396E-3 \& $2.65$ & 0.45323E-4 \& $3.16$ & 0.44637E-5 \& $3.34$ & 0.13580E-5 \& $1.72$
\tabularnewline &\\ [-.3cm] \hline &\\ [-.3cm]
$\parallel \bar{e}_h \parallel_{0,h}$ \& $\bar{r}_{0}$ & $\rightarrow$ & 0.27043E-1 \& $1.74$ & 0.70830E-2 \& $1.93$ & 0.17919E-2 \& $1.98$ & 0.44939E-3 \& $2.00$ 
\tabularnewline &\\ [-.3cm] \hline &\\ [-.3cm] 
\tabularnewline &\\ [-.3cm] \hline &\\ [-.3cm]
$| e_h |_{1,h}$ \& $r_{1}$& $\rightarrow$ & 0.39119E-6 \& $--$ & 0.13594E-5 \& $--$ & 0.94788E-5 \& $--$ & 0.21501E-4 \& $--$  
\tabularnewline &\\ [-.3cm] \hline &\\ [-.3cm] 
$| \tilde{e}_{h}|_{1,\tilde{h}}$ \& $\tilde{r}_1$ & $\rightarrow$ & 0.99920E-2 \& $1.93$ & 0.20221E-2 \& $2.30$ & 0.38103E-3 \& $2.41$ & 0.70845E-4 \& $2.43$ 
\tabularnewline &\\ [-.3cm] \hline &\\ [-.3cm] 
$| \bar{e}_{h} |_{1,h}$ \& $\bar{r}_1$ & $\rightarrow$ & 0.21044E-1 \& $1.60$ & 0.56045E-2 \& $1.91$ & 0.14255E-2 \& $1.98$ & 0.35878E-3 \& $1.99$ 
\tabularnewline &\\ [-.3cm] \hline &\\ [-.3cm]
\tabularnewline &\\ [-.3cm] \hline &\\ [-.3cm] 
$| e_h |_{\infty,h}$ \& $r_{\infty}$ & $\rightarrow$ & 0.87463E-6 \& $--$ & 0.38417E-5 \& $--$ & 0.19209E-4 \& $--$ & 0.29365E-4 \& $--$ 
\tabularnewline &\\ [-.3cm] \hline &\\ [-.3cm] 
$| \tilde{e}_{h} |_{\infty,\tilde{h}}$ \& $\tilde{r}_{\infty}$ & $\rightarrow$ & 0.10359E-1 \& $2.43$ & 0.14411E-2 \& $2.85$ & 0.17941E-3 \& $3.01$ & 0.25712E-4 \& $2.80$
\tabularnewline &\\ [-.3cm] \hline &\\ [-.3cm] 
$| \bar{e}_{h} |_{\infty,h}$ \& $\bar{r}_{\infty}$ & $\rightarrow$ & 0.11044E+0\& $1.79$ & 0.27032E-1 \& $2.03$ & 0.65678E-2 \& $2.04$ & 0.16214E-2 \& $2.02$ 
\tabularnewline &\\ [-.3cm] \hline &\\ [-.3cm] 
\end{tabular}
%\vskip.2cm 
\caption{Errors \& $IOC$'s for Test-problem 2 solved by \eqref{vph}, isoparametrics and \eqref{barvph} with $k=2$} 
}
\label{table2}
\end{table}

\item 

\underline{\bf Test-problem 3:} $\Omega=\Omega_{\beta}$ with $\beta=0.25$, $r \equiv 1$, $\gamma \equiv 1$ and $u=\rho^2$.\\

Akin to the previous test-problem, in this case we checked the consistency of the three methods in the presence of Robin boundary conditions, 
but now in the case of a non convex domain. From the results displayed in Table 3 we can assert that all the three methods show a regular behavior all the way. In the case of isoparametrics this is in contrast to Test-problems 1 and 2, since the $IOC$ downgrade is milder for the finest mesh. This effect could be due to the fact that the corresponing absolute errors are significantly larger here. Indeed, it might be that the degree of magnitude for which such a deterioration takes place - eventually owing to the stop criterion of the solution algorithm -, was not attained here for the finest mesh.
\begin{table}[h!]
{\small 
\centering
\begin{tabular}{cccccc} &\\ [-.3cm]  
$h$ & $\rightarrow$ & $1/2$ & $1/4$ & $1/8$ & $1/16$  
\tabularnewline &\\ [-.3cm] \hline &\\ [-.3cm]
$\parallel e_h \parallel_{0,h}$ \& $r_0$ & $\rightarrow$ &  0.16780E-6 \& $--$ &  0.38373E-6 \& $--$ &  0.82852E-6 \& $--$ &  0.58123E-5 \& $--$
\tabularnewline &\\ [-.3cm] \hline &\\ [-.3cm] 
$\parallel \tilde{e}_h \parallel_{0,\tilde{h}}$ \& $\tilde{r}_0$ & $\rightarrow$ &  0.17042E-2 \& $2.11$ &  0.12207E-3 \& $3.80$ &  0.11165E-4 \& $3.45$ &  0.13700E-5 \& $3.03$
\tabularnewline &\\ [-.3cm] \hline &\\ [-.3cm] 
$\parallel \bar{e}_h \parallel_{0,h}$ \& $\bar{r}_{0}$ & $\rightarrow$ &  0.38499E-1 \& $2.02$ &  0.98408E-2 \& $1.97$ &  0.24948E-2 \& $1.98$ &  0.62606E-3 \& $1.99$ 
\tabularnewline &\\ [-.3cm] \hline &\\ [-.3cm] 
\tabularnewline &\\ [-.3cm] \hline &\\ [-.3cm]
$| e_h |_{1,h}$ \& $r_{1}$& $\rightarrow$ &  0.37128E-6 \& $--$ &  0.86533E-6 \& $--$ &  0.42989E-5 \& $--$ &  0.13207E-4 \& $--$ 
\tabularnewline &\\ [-.3cm] \hline &\\ [-.3cm] 
$| \tilde{e}_{h}|_{1,\tilde{h}}$ \& $\tilde{r}_1$ & $\rightarrow$ &  0.19764E-1 \& $1.77$ &  0.40863E-2 \& $2.27$ &  0.76459E-3 \& $2.42$ &  0.13963E-3 \& $2.45$ 
\tabularnewline &\\ [-.3cm] \hline &\\ [-.3cm] 
$| \bar{e}_{h} |_{1,h}$ \& $\bar{r}_1$ & $\rightarrow$ &  0.52942E-1 \& $1.20$ &  0.15413E-1 \& $1.78$ &  0.39343E-2 \& $1.97$ &  0.98618E-3 \& $2.00$ 
\tabularnewline &\\ [-.3cm] \hline &\\ [-.3cm]
\tabularnewline &\\ [-.3cm] \hline &\\ [-.3cm] 
$| e_h |_{\infty,h}$ \& $r_{\infty}$ & $\rightarrow$ &  0.84090E-6 \& $--$ &  0.11345E-5 \& $-.43$ &  0.64039E-5 \& $--$ &  0.14681E-4 \& $--$ 
\tabularnewline &\\ [-.3cm] \hline &\\ [-.3cm] 
$| \tilde{e}_{h} |_{\infty,\tilde{h}}$ \& $\tilde{r}_{\infty}$ & $\rightarrow$ &  0.15588E-1 \& $1.73$ &  0.22274E-2 \& $2.81$ &  0.27562E-3 \& $3.01$ &  0.29344E-4 \& $3.23$
\tabularnewline &\\ [-.3cm] \hline &\\ [-.3cm] 
$| \bar{e}_{h} |_{\infty,h}$ \& $\bar{r}_{\infty}$ & $\rightarrow$ &  0.13219E+0\& $1.64$ &  0.31310E-1 \& $2.08$ &  0.74505E-2 \& $2.07$ &  0.18072E-2 \& $2.04$ 
\tabularnewline &\\ [-.3cm] \hline &\\ [-.3cm] 
\end{tabular}
%\vskip.2cm 
\caption{Errors \& $IOC$'s for Test-problem 3 solved by \eqref{vph}, isoparametrics and \eqref{barvph} with $k=2$} 
}
\label{table3}
\end{table}

\end{itemize}
It is noteworthy that the isoparametric method showed up significally more accurate than the do-nothing strategy in all the test-cases above.\\

We pursue this method assessment with test-problems having a non quadratic exact solution. \\
For convenience we define the quantity 
$$\eta = \sqrt{x^2/a^2+y^2/b^2+z^2}.$$

\begin{itemize}
\item 

\underline{\bf Test-problem 4:} $\Omega=\Omega_{ab}$ with $a=0.6$ and $b=0.8$, $r \equiv 1$, $\gamma \equiv 0$ and $u=\eta^2/2-\eta^3/3$.\\

Here we assess the three methods in the case of a non polynomial exact solution satisfying homogeneous Neumann boundary conditions. As we can infer from Table 4, the expected orders of convergence are reflected by the evolution of the $IOC$, including the one of the method \eqref{barvph} in terms of the $| \cdot |_{1,h}$. Indeed, from an argument developed in \cite{Arxiv2026} the corresponding order of convergence is supposed to be two in this case, since the exact solution is constant on the boundary. But what is really surprising is the accuracy of the do-nothing approach under this criterion, since it is not so different from the one the two other methods, even if slightly worse, while a small advantage prevails for $e_h$ over $\tilde{e}_h$. On the other hand, regarding the two other criteria the errors for the isoparametric method leave no room for doubts about its superiority in terms of accuracy, in spite of a persistent $IOC$ downgrade as the mesh is refined.  %.                                                                                                                                       

\begin{table}[h!]
{\small 
\centering
\begin{tabular}{cccccc} &\\ [-.3cm]  
$h$ & $\rightarrow$ & $1/2$ & $1/4$ & $1/8$ & $1/16$  
\tabularnewline &\\ [-.3cm] \hline &\\ [-.3cm]
$\parallel e_h \parallel_{0,h}$ \& $r_0$ & $\rightarrow$ &  0.28684E-1 \& $2.39$ &  0.40409E-2 \& $2.83$ &  0.55514E-3 \& $2.86$ &  0.72222E-4 \& $2.94$
\tabularnewline &\\ [-.3cm] \hline &\\ [-.3cm] 
$\parallel \tilde{e}_h \parallel_{0,\tilde{h}}$ \& $\tilde{r}_0$ & $\rightarrow$ &  0.66233E-03 \& $1.61$ &  0.96397E-04 \& $2.78$ &  0.12632E-04 \& $2.93$ &  0.21030E-05 \& $2.59$
\tabularnewline &\\ [-.3cm] \hline &\\ [-.3cm] 
$\parallel \bar{e}_h \parallel_{0,h}$ \& $\bar{r}_{0}$ & $\rightarrow$ &  0.57427E-1 \& $1.72$ &  0.15151E-1 \& $1.92$ &  0.38400E-2 \& $1.98$ &  0.96333E-3 \& $2.00$ 
\tabularnewline &\\ [-.3cm] \hline &\\ [-.3cm] 
\tabularnewline &\\ [-.3cm] \hline &\\ [-.3cm]
$| e_h |_{1,h}$ \& $r_{1}$& $\rightarrow$ &  0.13910E-1 \& $1.48$ &  0.38595E-2 \& $1.85$ &  0.10267E-2 \& $1.91$ &  0.26630E-3 \& $1.95$ 
\tabularnewline &\\ [-.3cm] \hline &\\ [-.3cm] 
$| \tilde{e}_{h}|_{1,\tilde{h}}$ \& $\tilde{r}_1$ & $\rightarrow$ &  0.14292E-1 \& $1.23$ &  0.40665E-2 \& $1.81$ &  0.10669E-2 \& $1.93$ &  0.27208E-3 \& $1.97$ 
\tabularnewline &\\ [-.3cm] \hline &\\ [-.3cm] 
$| \bar{e}_{h} |_{1,h}$ \& $\bar{r}_1$ & $\rightarrow$ &  0.17235E-1 \& $1.36$ &  0.49737E-2 \& $1.79$ &  0.13178E-2 \& $1.92$ &  0.33868E-3 \& $1.96$ 
\tabularnewline &\\ [-.3cm] \hline &\\ [-.3cm]
\tabularnewline &\\ [-.3cm] \hline &\\ [-.3cm] 
$| e_h |_{\infty,h}$ \& $r_{\infty}$ & $\rightarrow$ &  0.66922E-1 \& $2.47$ &  0.93247E-2 \& $2.84$ &  0.13091E-2 \& $2.83$ &  0.17122E-3 \& $2.93$ 
\tabularnewline &\\ [-.3cm] \hline &\\ [-.3cm] 
$| \tilde{e}_{h} |_{\infty,\tilde{h}}$ \& $\tilde{r}_{\infty}$ & $\rightarrow$ &  0.60343E-2 \& $3.46$ &  0.85220E-3 \& $2.82$ &  0.12474E-3 \& $2.77$ &  0.18592E-4 \& $2.75$
\tabularnewline &\\ [-.3cm] \hline &\\ [-.3cm] 
$| \bar{e}_{h} |_{\infty,h}$ \& $\bar{r}_{\infty}$ & $\rightarrow$ &  0.13128E+0 \& $1.82$ &  0.34119E-1 \& $1.94$ &  0.86114E-2 \& $1.99$ &  0.21528E-2 \& $2.00$ 
\tabularnewline &\\ [-.3cm] \hline &\\ [-.3cm] 
\end{tabular}
%\vskip.2cm 
\caption{Errors \& $IOC$'s for Test-problem 4 solved by \eqref{vph}, isoparametrics and \eqref{barvph} with $k=2$} 
}
\label{table4}
\end{table}

\item 

\underline{\bf Test-problem 5:} $\Omega=\Omega_{ab}$ with $a=0.6$ and $b=0.8$, $r \equiv 0$, $\gamma \equiv 1$ and $u=\eta^2/2-\eta^3/3$.\\

In this test-case we took the same domain and the same exact solution as in the previous one. However, here we solved a pure diffusion problem, more specifically, a Poisson equation with inhomogeneous boundary conditions. In spite of the problem's different nature, in regard of the results supplied in Table 5, there is no significant difference in the behavior of the three methods being compared with respect to Test-problem 4. It should be underlined, however, that the gradual deterioration of the $IOC$s for the isoparametric solution as the mesh is refined is less conspicuous here than in Test-problems 1,2 and 4, although the magnitude of the corresponding absolute errors is roughly the same in all cases.  
  
\begin{table}[h!]
{\small 
\centering
\begin{tabular}{cccccc} &\\ [-.3cm]  
$h$ & $\rightarrow$ & $1/2$ & $1/4$ & $1/8$ & $1/16$  
\tabularnewline &\\ [-.3cm] \hline &\\ [-.3cm]
$\parallel e_h \parallel_{0,h}$ \& $r_0$ & $\rightarrow$ &  0.74861E-2 \& $2.43$ &  0.10345E-2 \& $2.86$ &  0.14156E-3 \& $2.87$ &  0.18512E-4 \& $2.93$
\tabularnewline &\\ [-.3cm] \hline &\\ [-.3cm] 
$\parallel \tilde{e}_h \parallel_{0,\tilde{h}}$ \& $\tilde{r}_0$ & $\rightarrow$ &  0.66873E-3 \& $1.61$ &  0.96693E-4 \& $2.79$ &  0.12633E-4 \& $2.94$ &  0.18131E-5 \& $2.80$
\tabularnewline &\\ [-.3cm] \hline &\\ [-.3cm] 
$\parallel \bar{e}_h \parallel_{0,h}$ \& $\bar{r}_{0}$ & $\rightarrow$ &  0.13791E-1 \& $1.55$ &  0.37504E-2 \& $1.88$ &  0.95778E-3 \& $1.97$ &  0.24084E-3 \& $1.99$ 
\tabularnewline &\\ [-.3cm] \hline &\\ [-.3cm] 
\tabularnewline &\\ [-.3cm] \hline &\\ [-.3cm]
$| e_h |_{1,h}$ \& $r_{1}$& $\rightarrow$ &  0.12608E-1 \& $1.17$ &  0.37627E-2 \& $1.74$ &  0.10198E-2 \& $1.88$ &  0.26573E-3 \& $1.94$ 
\tabularnewline &\\ [-.3cm] \hline &\\ [-.3cm] 
$| \tilde{e}_{h}|_{1,\tilde{h}}$ \& $\tilde{r}_1$ & $\rightarrow$ &  0.14293E-1 \& $1.23$ &  0.40666E-2 \& $1.81$ &  0.10669E-2 \& $1.93$ &  0.27204E-3 \& $1.97$ 
\tabularnewline &\\ [-.3cm] \hline &\\ [-.3cm] 
$| \bar{e}_{h} |_{1,h}$ \& $\bar{r}_1$ & $\rightarrow$ &  0.12747E-1 \& $1.16$ &  0.38314E-2 \& $1.73$ &  0.10383E-2 \& $1.88$ &  0.27033E-3 \& $1.94$ 
\tabularnewline &\\ [-.3cm] \hline &\\ [-.3cm]
\tabularnewline &\\ [-.3cm] \hline &\\ [-.3cm] 
$| e_h |_{\infty,h}$ \& $r_{\infty}$ & $\rightarrow$ &  0.20748E-1 \& $2.70$ &  0.30001E-2 \& $2.79$ &  0.42905E-3 \& $2.81$ &  0.57290E-4 \& $2.90$ 
\tabularnewline &\\ [-.3cm] \hline &\\ [-.3cm] 
$| \tilde{e}_{h} |_{\infty,\tilde{h}}$ \& $\tilde{r}_{\infty}$ & $\rightarrow$ &  0.60151E-2 \& $3.46$ &  0.85264E-3 \& $2.82$ &  0.12053E-3 \& $2.82$ &  0.18121E-4 \& $2.73$
\tabularnewline &\\ [-.3cm] \hline &\\ [-.3cm] 
$| \bar{e}_{h} |_{\infty,h}$ \& $\bar{r}_{\infty}$ & $\rightarrow$ &  0.33638E-1 \& $2.01$ &  0.88645E-2 \& $1.92$ &  0.22414E-2 \& $1.98$ &  0.56296E-3 \& $1.99$ 
\tabularnewline &\\ [-.3cm] \hline &\\ [-.3cm] 
\end{tabular}
%\vskip.2cm 
\caption{Errors \& $IOC$'s for Test-problem 5 solved by \eqref{vph}, isoparametrics and \eqref{barvph} with $k=2$} 
}
\label{table5}
\end{table}

\item 

\underline{\bf Test-problem 6:} $\Omega=\Omega_{\beta}$ with $\beta=0.25$, $r \equiv 1$, $\gamma \equiv 1$ and $u=\rho^4/4.$.

The main purpose of this test-problem is to check the performance of the three methods under study, in the solution of a reaction-diffusion equation with (inhomogeneous) Robin boundary conditions in a non convex domain. We took a polynomial of degree four as an exact solution.\\
Regarding the measures $\| \cdot \|_{0,h}$ and $| \cdot |_{\infty,h}$, here again the isoparametric solution is more accurate than the Petrov-Galerkin one. However the discrepancy here are lesser than in the two previous test-cases, especially in the latter measure. On the other hand, no real $IOC$ downgrade for the finest meshes is observed here for isoparametrics, eventually because the magnitude of the corresponding errors is larger in this case. \\
It is also interesting to highlight the overall more stable evolution the of $IOC$s for the Petrov-Galerkin solution $u_h$, and also for $\bar{u}_h$ to some extent, in contrast to the more random one for the isoparametric solution $\tilde{u}_h$. This behavior can also be observed in the case of all the other test-problems addressed in this section, wherever applicable.
\begin{table}[h!]
{\small 
\centering
\begin{tabular}{cccccc} &\\ [-.3cm]  
$h$ & $\rightarrow$ & $1/2$ & $1/4$ & $1/8$ & $1/16$  
\tabularnewline &\\ [-.3cm] \hline &\\ [-.3cm]
$\parallel e_h \parallel_{0,h}$ \& $r_0$ & $\rightarrow$  &  0.15525E-1 \& $2.21$ &  0.23743E-2 \& $2.71$ &  0.32104E-3 \& $2.89$ &  0.45136E-4 \& $2.83$
\tabularnewline &\\ [-.3cm] \hline &\\ [-.3cm] 
$\parallel \tilde{e}_h \parallel_{0,\tilde{h}}$ \& $\tilde{r}_0$ & $\rightarrow$ &  0.28655E-2 \& $1.60$ &  0.39447E-3 \& $2.86$ &  0.48370E-4 \& $3.03$ &  0.61431E-5 \& $2.98$
\tabularnewline &\\ [-.3cm] \hline &\\ [-.3cm] 
$\parallel \bar{e}_h \parallel_{0,h}$ \& $\bar{r}_{0}$ & $\rightarrow$ &  0.57427E-1 \& $1.72$ &  0.15151E-1 \& $1.92$ &  0.38400E-2 \& $1.98$ &  0.96333E-3 \& $2.00$ 
\tabularnewline &\\ [-.3cm] \hline &\\ [-.3cm] 
\tabularnewline &\\ [-.3cm] \hline &\\ [-.3cm]
$| e_h |_{1,h}$ \& $r_{1}$& $\rightarrow$ &  0.43719E-1 \& $1.72$ &  0.11811E-1 \& $1.89$ &  0.30207E-2 \& $1.97$ &  0.75955E-3 \& $1.99$ 
\tabularnewline &\\ [-.3cm] \hline &\\ [-.3cm] 
$| \tilde{e}_{h}|_{1,\tilde{h}}$ \& $\tilde{r}_1$ & $\rightarrow$ &  0.32073E-1 \& $0.98$ &  0.94046E-2 \& $1.77$ &  0.25453E-2 \& $1.89$ &  0.66183E-3 \& $1.94$ 
\tabularnewline &\\ [-.3cm] \hline &\\ [-.3cm] 
$| \bar{e}_{h} |_{1,h}$ \& $\bar{r}_1$ & $\rightarrow$ &  0.57018E-1 \& $0.99$ &  0.16803E-1 \& $1.76$ &  0.44176E-2 \& $1.93$ &  0.11242E-2 \& $1.97$ 
\tabularnewline &\\ [-.3cm] \hline &\\ [-.3cm]
\tabularnewline &\\ [-.3cm] \hline &\\ [-.3cm] 
$| e_h |_{\infty,h}$ \& $r_{\infty}$ & $\rightarrow$ &  0.31892E-1 \& $2.46$ &  0.52978E-2 \& $2.59$ &  0.79033E-3 \& $2.74$ &  0.11537E-3 \& $2.78$ 
\tabularnewline &\\ [-.3cm] \hline &\\ [-.3cm] 
$| \tilde{e}_{h} |_{\infty,\tilde{h}}$ \& $\tilde{r}_{\infty}$ & $\rightarrow$ &  0.16399E-1 \& $3.03$ &  0.32849E-2 \& $2.32$ &  0.46151E-3 \& $2.83$ &  0.65794E-4 \& $2.81$
\tabularnewline &\\ [-.3cm] \hline &\\ [-.3cm] 
$| \bar{e}_{h} |_{\infty,h}$ \& $\bar{r}_{\infty}$ & $\rightarrow$ &  0.11820E+0 \& $1.48$ &  0.31962E-1 \& $1.89$ &  0.79846E-2 \& $2.00$ &  0.19908E-2 \& $2.00$ 
\tabularnewline &\\ [-.3cm] \hline &\\ [-.3cm] 
\end{tabular}
%\vskip.2cm 
\caption{Errors \& $IOC$'s for Test-problem 6 solved by \eqref{vph}, isoparametrics and \eqref{barvph} with $k=2$} 
}
\label{table6}
\end{table}

\end{itemize}

\subsection{Other comparison criteria}

Accuracy and order of convergence are certainly the most important focus in comparing the performances of the different approximation methods in this study. 
Nevertheless it is instructive to say a word on complexity and processing time, since the effective solution procedure is quite different for one them. As a mattter of fact, it might seem that the new Petrov-Galerkin formulation is more complex and more time consuming than the two others. Actually, although its implementation is a little harder, this does not really induce any significant increase of running time, at least as long as the iterative solution procedure \eqref{iterproc} combined with the conjugate gradient method is employed. \\
In this respect, first of all it is important to recall the argument in \cite{JCAMbis}, according to which, the finer the mesh is the faster this procedure converges. For example, while for the very coarse mesh with only 6 tetrahedra seven iterations are necessary to satisfy the tolerance of $10^{-7}$ for the measure $| \cdot |_{\infty,h}$, this number decreases gradually as the mesh is refined, to become only three for the finest mesh with almost 25,000 tetrahedra.  Furthermore the number of iterations necessary to satisfy the stop criterion for the conjugate gradient method may vary in an unpredictable manner for each method, which can lead to unexpected results in terms of processing  time. Just to illustrate this kind of effect, we supply in Table 7 running times in seconds to solve one of the above test-problems as representative of all the others. It is particularly noteworthy that our Petrov-Galerkin approach is not at all more time-consuming than the 
isoparametric method and the do-nothing approach is not cheaper than the two others in terms of running time all the way.      

\begin{table*}[h!]
{\small 
\centering
\begin{tabular}{cccccc} &\\ [-.3cm]  
$h$ & $\longrightarrow$ & $1/2$ & $1/4$ & $1/8$ & $1/16$  
\tabularnewline &\\ [-.3cm] \hline &\\ [-.3cm]
Petrov-Galerkin method (13) & $\longrightarrow$ & 0.1470 & 1.0189 & 20.0186 & 908.2824  
\tabularnewline &\\ [-.3cm] \hline &\\ [-.3cm] 
Isoparamtric method         & $\longrightarrow$ & 0.1253 & 1.0208 & 22.4365 & 904.4238 
\tabularnewline &\\ [-.3cm] \hline &\\ [-.3cm] 
Galerkin method (18)        & $\longrightarrow$ & 0.1098 & 0.8860 & 19.7942 & 939.8917
\tabularnewline &\\ [-.3cm] \hline &\\ [-.3cm] 
\end{tabular}
%\vskip.2cm 
\caption{Typical running times in seconds to solve a test-problem with the three assessed methods} 
}
\label{tabletime}
\end{table*}

\section{Conclusions and perspectives}
   
As a brief summary of this work we list below some conclusions, mostly drawn on the basis of the results supplied in Section 4.
\begin{enumerate}
\item The new Petrov-Galerkin formulation defined in Section 3 in connection with Lagrange tetrahedral finite elements of degree $k >1$      for second order boundary value problems in smooth domains with natural conditions, showed to be a valid alternative to the well-established isoparametric technique of the same degree, at least for $k=2$. Indeed, they yield approximations roughly of the same order in the three most common norm and semi-norms at a comparable cost. Of course, this was not the case of the do-nothing strategy \eqref{barvph}, although, surprisingly enough, this approach yielded roughly second order approximations in the $| \cdot |_{1,h}$ semi-norm, even in the cases where the exact solution's tangential derivatives do not vanish everywhere on $\Gamma$.
\vspace{-1mm}
\item The new method appeared to be more stable than the isoparametric method in terms of the variation of the convergence rates as the meshes are refined. Anyway, in many cases some erosion was observed for both methods mostly with respect to the expected third order convergence in the measures $\| \cdot \|_{0,h}$ and $| \cdot |_{\infty,h}$. The sub-optimal second order Galerkin method \eqref{barvph} in turn behaved fairly well in this respect.
\vspace{-1mm}  
\item
While isoparametrics carry the intrinsic complexity of rational shape functions, the Petrov-Galerkin formulation \eqref{vph} benefits from the simplicity inherent to polynomial algebra. In particular this allows the user of the latter method to do without numerical quadrature to approximate the terms of element matrices and vectors consisting of integrals in the tetrahedra. The burden of selecting a suitable quadrature formula for these computations can be thus left to integrals on the plane triangular boundary faces. However, here again, this can be replaced by a classical Lagrange interpolation of the problem's data, such as the direction cosines of the outer normal vector to the curvilinear boundary or the boundary itself.           
\vspace{-1mm}                                    
\item 
The isoparametric solutions turned out to be fairly more accurate than the Petrov-Galerkin's, as far as the two measures not involving gradients are concerned. On the other hand both methods appeared to be equally accurate in the semi-norm $| \cdot |_{1,h}$, with a slight advantage of the latter.
\end{enumerate}

In complement to the above discussion, the author underlines that the confrontation of the classical isoparametric method with the Petrov-Galerkin approach based on piecewise polynomial shape-functions upon affine fitted simplicial meshes cannot be seen as settled is this sole article. Neither in \cite{Arxiv2017}, \cite{ZAMM}, \cite{IMAJNA} and \cite{Arxiv2026}, in which it was also carried out for second order methods in the energy norm in the case of second order boundary value problems in two- or three-dimensional domains with both essential and natural boundary conditions. In this precise framework, the experiments indicated that both approaches can be viewed as globally equivalent. Indeed, depending on the case or the criterion, one of them showed up superior to the other. \\
The assessment of both techniques could be pushed further, for instance, by comparing them in association with higher order interpolations. However, from the accumulated experience it is clear beforehand that the well-established isoparametric finite element method is quite effective for the aforementioned class of partial differential equations. Nevertheless, the rather recent technique studied in the four papers quoted above can be viewed as an attractive alternative to practitioners who wish to take advantage of some of its features already highlighted in this article, as compared to the isoparametric method. \\
Incidentally, it should be pointed out that one of the main merits of the former approach as compared to the isoparametric one is its universality. This is because picking up the prescribed conditions on the true curved boundary, while using a fitted straight-edged mesh is very intuitive and hence can be quite easily applied to finite element methods for more challenging types of equations. In this respect  fourth order elliptic equations in two-dimension space such as those governing thin plate bending are an outstanding example. Indeed, in this case our approach circumvents the need for ad hoc techniques such as the one in \cite{Mansfield} for the so-called Clough-Tocher triangle, since known parametric counterparts of methods for polygonal domains are limited to a certain type of finite elements (see e.g. \cite{Zenisek}). For this reason the author intends to extend in this direction his studies on the newer technique in the near future.  

%universality of the principle the technique adopted in this article is based upon, to handle essential conditions prescribed on curvilinear boundaries. In this sense, its application is not restricted at all to finite element methods. Actually, successful adaptions thereof to other numerical methods commonly used in Engineering and Science, such as the discontinuous Galerkin method, can be found in the literature. %(see e.g. \cite{DG2D} and references therein).      

%\section*{Acknowledgement} 
%\vspace*{-3mm} 
%The author is thankful to Fleurianne Bertrand for helpful discussions.

%\section*{Declarations}
%\begin{itemize}
%\item
%Data sets including codes generated during the current study are available from the author on reasonable request. %These include the Fortran 95 codes developed to the their generation itself.  
%\item
%The author has no relevant financial or non-financial interests to disclose and has received no funding for conducting this study.
%\end{itemize}

\end{document}